# STEPUP PROCEDURES FOR CONTROL OF GENERALIZATIONS OF THE FAMILYWISE ERROR RATE

BY JOSEPH P. ROMANO AND AZEEM M. SHAIKH

*Stanford University*

Consider the multiple testing problem of testing null hypotheses $H_1, \ldots, H_s$. A classical approach to dealing with the multiplicity problem is to restrict attention to procedures that control the familywise error rate (*FWER*), the probability of even one false rejection. But if $s$ is large, control of the *FWER* is so stringent that the ability of a procedure that controls the *FWER* to detect false null hypotheses is limited. It is therefore desirable to consider other measures of error control. This article considers two generalizations of the *FWER*. The first is the $k$-*FWER*, in which one is willing to tolerate $k$ or more false rejections for some fixed $k \geq 1$. The second is based on the false discovery proportion (*FDP*), defined to be the number of false rejections divided by the total number of rejections (and defined to be 0 if there are no rejections). Benjamini and Hochberg [*J. Roy. Statist. Soc. Ser. B* **57** (1995) 289–300] proposed control of the false discovery rate (*FDR*), by which they meant that, for fixed $\alpha$, $E(FDP) \leq \alpha$. Here, we consider control of the *FDP* in the sense that, for fixed $\gamma$ and $\alpha$, $P\{FDP > \gamma\} \leq \alpha$. Beginning with any nondecreasing sequence of constants and $p$-values for the individual tests, we derive stepup procedures that control each of these two measures of error control without imposing any assumptions on the dependence structure of the $p$-values. We use our results to point out a few interesting connections with some closely related stepdown procedures. We then compare and contrast two *FDP*-controlling procedures obtained using our results with the stepup procedure for control of the *FDR* of Benjamini and Yekutieli [*Ann. Statist.* **29** (2001) 1165–1188].

**1. Introduction.** In this article we consider the problem of simultaneously testing hypotheses $H_i$ ($i = 1, \ldots, s$). We shall assume that tests based on $p$-values $\hat{p}_1, \ldots, \hat{p}_s$ are available for the individual hypotheses and that the question of interest is how to combine these $p$-values into a simultaneous testing procedure. In other words, each $\hat{p}_i$ is a marginal $p$-value in the sense that it could be used for testing $H_i$; $p$-values for testing individual hypotheses are reviewed in Section 3.3 of [11].

A classical approach to handling the multiplicity problem is to restrict attention to procedures that control the *familywise error rate* (*FWER*), defined to be the probability of one or more false rejections. When evaluating a testing procedure, one must consider not only control of false rejections, but also the ability of the procedure to detect departures from the null hypothesis when they do occur. When the number of tests $s$ is large, control of the *FWER* is so stringent that departures from the null hypothesis have little chance of being detected. As a result, alternative measures of error control have been considered, which control false rejections less severely, but in doing so are better able to detect false null hypotheses.

Hommel and Hoffman [8] and Lehmann and Romano [10] considered control of the $k$-*FWER*, the probability of rejecting at least $k$ true null hypotheses. Such an error rate with $k > 1$ is appropriate when one is willing to tolerate one or more false rejections, provided the number of false rejections is controlled. Evidently, taking $k = 1$ reduces to the usual *FWER*. These authors derived both single step and stepdown methods that guarantee that the $k$-*FWER* is bounded above by $\alpha$.

Lehmann and Romano [10] also considered control of the *false discovery proportion* (*FDP*), defined as the total number of false rejections divided by the total number of rejections (and equal to 0 if there are no rejections). Given a user specified value $\gamma \in [0, 1]$, control of the *FDP* means we wish to ensure that $P\{FDP > \gamma\}$ is bounded above by $\alpha$. Setting $\gamma = 0$ reduces to the usual *FWER*. Lehmann and Romano [10] also provided stepdown procedures for control of the *FDP* that hold under either mild or no assumptions on the joint distribution of the $p$-values. Romano and Shaikh [12] improved upon these arguments to derive a stepdown procedure for control of the *FDP* that is also valid under no restrictions on the dependence structure of the $p$-values, but considerably more powerful than the method proposed in [10]. In this article, unlike either of these previous works, we consider stepup procedures. We derive stepup procedures that control the $k$-*FWER* and the *FDP* under no assumptions on the joint distribution of the $p$-values.

A closely related type of error control that has received much attention since it was first proposed in [1] is control of the *false discovery rate* (*FDR*), which demands that $E(FDP)$ is bounded above by $\alpha$. This original paper imposed the very strong assumption that the $p$-values were independent, but Benjamini and Yekutieli [2] have since proposed a stepup method that is valid under no assumptions on the joint distribution of the $p$-values. It is of interest to compare control of the *FDP* with control of the *FDR*. Even though ensuring that the *FDR* is bounded does not prohibit the *FDP* from varying, some obvious connections between methods that control the *FDP* in the sense that

$$P\{FDP > \gamma\} \leq \alpha \tag{1}$$

and methods the control its expected value, the *FDR*, can be made. Indeed, for any random variable $X$ on $[0,1]$, we have

$$E(X) = E(X|X \leq \gamma)P\{X \leq \gamma\} + E(X|X > \gamma)P\{X > \gamma\}$$
$$\leq \gamma P\{X \leq \gamma\} + P\{X > \gamma\},$$

which leads to

$$\frac{E(X) - \gamma}{1 - \gamma} \leq P\{X > \gamma\} \leq \frac{E(X)}{\gamma}, \tag{2}$$

with the last inequality just Markov's inequality. Applying this to $X = FDP$, we see that, if a method controls the *FDR* at level $q$, then it controls the *FDP* in the sense $P\{FDP > \gamma\} \leq q/\gamma$. Conversely, if the *FDP* is controlled in the sense of (1), then the *FDR* is controlled at level $\gamma(1-\alpha) + \alpha$. Therefore, in principle, a method that controls the *FDP* in the sense of (1) can be used to control the *FDR* and vice versa, as previously noted by van der Laan, Dudoit and Pollard [16]. We will compare methods for control of the *FDP* with the method for control of the *FDR* proposed by Benjamni and Yekutieli [2] in light of this observation. Note that setting $\alpha = 1/2$ restricts the median of the *FDP* to be no greater than $\gamma$.

A growing literature has proposed various procedures which control generalized error rates. Genovese and Wasserman [4], for example, study asymptotic procedures that control the *FDP* and the *FDR* in the framework of a random effects mixture model. These ideas are extended in [4]. Korn, Troendle, McShane and Simon [9] provide methods that control both the *k-FWER* and *FDP*; their results are limited to a multivariate permutation model. Their results are generalized in [13]. Alternative procedures for control of the *k-FWER* and *FDP* are given in [16].

The paper is organized as follows. In Section 2 we describe our terminology and the class of stepup procedures. All of our methods assume that marginal $p$-values are available for testing each of the individual hypotheses, in the sense described in (3). Our methods are designed to hold under no dependence assumptions among the $p$-values, but do not attempt to estimate the dependence structure (as in van der Laan, Dudoit and Pollard [16] or Romano and Wolf [13]). Hence, our main results are exact and nonasymptotic; however, if the individual $p$-values are only approximate (as they typically are when using asymptotic approximations or resampling methods), the error control will hold approximately; see Remark 4.2. Control of the *k-FWER* and *FDP* are considered, respectively, in Sections 3 and 4. Our calculations in these two sections shed some light on the relationship between stepup and stepdown procedures as well. In Section 5 we use the relationship (2) to compare methods for controlling the *FDP* with the method of Benjamini and Yekutieli [2] for controlling the *FDR*. Section 6 illustrates the method with two examples. Section 7 concludes.

**2. A class of stepup procedures.** A formal description of our setup is as follows. Suppose data $X$ is available from some model $P \in \Omega$. A general hypothesis $H$ can be viewed as a subset $\omega$ of $\Omega$. For testing $H_i : P \in \omega_i$, $i = 1, \ldots, s$, let $I(P)$ denote the set of true null hypotheses when $P$ is the true probability distribution; that is, $i \in I(P)$ if and only if $P \in \omega_i$.

We assume that $p$-values $\hat{p}_1, \ldots, \hat{p}_s$ are available for testing $H_1, \ldots, H_s$. Specifically, we mean that $\hat{p}_i$ must satisfy

(3)     $P\{\hat{p}_i \leq u\} \leq u$     for any $u \in (0,1)$ and any $P \in \omega_i$.

Note that we do not require $\hat{p}_i$ to be uniformly distributed on $(0,1)$ if $H_i$ is true, in order to accommodate discrete situations. In deriving our results, we assume that (3) holds exactly, but we show in Remark 4.2 below that all of our results also extend to the case in which the $p$-values only satisfy (3) approximately.

In general, a $p$-value $\hat{p}_i$ will satisfy (3) if it is obtained from a nested set of rejection regions. In other words, suppose $S_i(\alpha)$ is a rejection region for testing $H_i$; that is,

(4)     $P\{X \in S_i(\alpha)\} \leq \alpha$     for all $0 < \alpha < 1, P \in \omega_i$

and $S_i(\alpha) \subset S_i(\alpha')$ whenever $\alpha < \alpha'$. Then the $p$-value $\hat{p}_i$ defined by $\hat{p}_i = \hat{p}_i(X) = \inf\{\alpha : X \in S_i(\alpha)\}$ satisfies (3). Such a construction applies to many parametric procedures and also some nonparametric procedures, such as those based on permutation or randomization tests; see (15.5) in [11].

In this article we consider the following class of *stepup* procedures. Let

(5)     $\alpha_1 \leq \alpha_2 \leq \cdots \leq \alpha_s$

be a nondecreasing sequence of constants. Order the $p$-values as

$$\hat{p}_{(1)} \leq \hat{p}_{(2)} \leq \cdots \leq \hat{p}_{(s)},$$

and let $H_{(1)}, \ldots, H_{(s)}$ denote the corresponding null hypotheses. If $\hat{p}_{(s)} \leq \alpha_s$, then reject all null hypotheses; otherwise, reject hypotheses $H_{(1)}, \ldots, H_{(r)}$, where $r$ is the smallest index satisfying

(6)     $\hat{p}_{(s)} > \alpha_s, \ldots, \hat{p}_{(r+1)} > \alpha_{r+1}.$

If, for all $r$, $\hat{p}_{(r)} > \alpha_r$, then reject no hypotheses. That is, a stepup procedure begins with the least significant $p$-value and continues accepting hypotheses as long as their corresponding $p$-values are large.

We will compare these stepup procedures considered with certain *stepdown* procedures. Given constants of the form (5), a stepdown procedure determines which null hypotheses to reject as follows. If $\hat{p}_{(1)} > \alpha_1$, then reject no null hypotheses; otherwise, reject hypotheses $H_{(1)}, \ldots, H_{(r)}$, where $r$ is the largest index satisfying

(7) $$\hat{p}_{(1)} \leq \alpha_1, \ldots, \hat{p}_{(r)} \leq \alpha_r.$$

That is, a stepdown procedure begins with the most significant $p$-value and continues rejecting hypotheses as long as their corresponding $p$-values are small.

REMARK 2.1. Consider a stepup and a stepdown procedure based on the same set of critical values (5). The stepup procedure will always reject at least as many hypotheses as the stepdown procedure. If both methods satisfy the given measure of error control, then the stepup procedure is more powerful than the corresponding stepdown procedure based on the same critical values in the sense that the stepup procedure will have a greater chance of detecting false null hypotheses.

**3. Control of the $k$-FWER.** In this section we consider control of the $k$-FWER, defined formally as

(8) $$P\{\text{reject} \geq k \text{ hypotheses } H_i \text{ with } i \in I(P)\}.$$

Control of the $k$-FWER at level $\alpha$ requires that $k$-FWER $\leq \alpha$ for all $P$. We first establish a result that will aid in constructing stepup methods that control the $k$-FWER.

LEMMA 3.1. *Consider testing $s$ null hypotheses, with $|I|$ of them true. Let*

$$\hat{q}_{(1)} \leq \cdots \leq \hat{q}_{(|I|)}$$

*denote the ordered values of the p-values corresponding to true hypotheses. Then, the stepup procedure based on constants $\alpha_1 \leq \cdots \leq \alpha_s$ satisfies*

(9) $$k - FWER \leq P\left\{ \bigcup_{k \leq j \leq |I|} \{\hat{q}_{(j)} \leq \alpha_{s-|I|+j}\} \right\}.$$

PROOF. Assume that $|I| \geq k$, for otherwise there is nothing to prove. Let $\hat{p}_{(1)} \leq \cdots \leq \hat{p}_{(s)}$ denote the ordered values of the $p$-values. For $1 \leq j \leq s$, let $A_j$ denote the event in which exactly $j$ hypotheses are rejected by the stepup procedure; that is,

$$A_j = \{\hat{p}_{(s)} > \alpha_s, \ldots, \hat{p}_{(j+1)} > \alpha_j, \hat{p}_{(j)} \leq \alpha_j\}.$$

Denote by $T$ the event in which at least $k$ true hypotheses are rejected. Consider the event $A_s \cap T$. Note that $A_s \cap T \subseteq \{\hat{p}_{(s)} \leq \alpha_s\} \cap T \subseteq \{\hat{q}_{(|I|)} \leq \alpha_s\}$. Likewise, note that $A_{s-1} \cap T \subseteq \{\hat{q}_{(|I|-1)} \leq \alpha_{s-1}\}$ if $|I| - 1 > k$ and $\subseteq \{\hat{q}_{(k)} \leq \alpha_{s-1}\}$ if $|I| - 1 \leq k$. In general, we have that

$$A_j \cap T \subseteq \begin{cases} \{\hat{q}_{(j+|I|-s)} \leq \alpha_j\}, & \text{if } j > s - |I| + k, \\ \{\hat{q}_{(k)} \leq \alpha_j\}, & \text{if } j \leq s - |I| + k. \end{cases}$$

Thus, the $k$-FWER is bounded above by the probability of the event

$$\bigcup_{k \leq j \leq s} A_j \cap T \subseteq \left\{ \bigcup_{k \leq j \leq s-|I|+k} \{\hat{q}_{(k)} \leq \alpha_j\} \right\} \cup \left\{ \bigcup_{s-|I|+k < j \leq s} \{\hat{q}_{(j+|I|-s)} \leq \alpha_j\} \right\}$$

$$\subseteq \bigcup_{s-|I|+k \leq j \leq s} \{\hat{q}_{(j+|I|-s)} \leq \alpha_j\} \subseteq \bigcup_{k \leq j \leq |I|} \{\hat{q}_{(j)} \leq \alpha_{s-|I|+j}\},$$

where the second inclusion follows from the fact that $\{\hat{q}_{(k)} \leq \alpha_j\} \subseteq \{\hat{q}_{(k)} \leq \alpha_{s-|I|+k}\}$ for $j \leq s - |I| + k$. The asserted claim now follows. $\square$

Given a sequence of constants $\alpha_1 \leq \cdots \leq \alpha_s$, we will now use Lemma 3.1 to construct a stepup procedure that controls the $k$-FWER. To this end, define

$$(10) \quad S_1 = S_1(k, s, |I|) = |I| \frac{\alpha_{s-|I|+k}}{k} + |I| \sum_{k < j \leq |I|} \frac{\alpha_{s-|I|+j} - \alpha_{s-|I|+j-1}}{j}$$

and let

$$(11) \quad D_1 = D_1(k, s) = \max_{k \leq |I| \leq s} S_1(k, s, |I|).$$

THEOREM 3.1. *Let $\alpha_1 \leq \cdots \leq \alpha_s$ be given. For testing $H_i : P \in \omega_i$, $i = 1, \ldots, s$, suppose $\hat{p}_i$ satisfies (3). Consider the stepup procedure with critical values $\alpha'_i = \alpha \alpha_i / D_1(k, s)$, where $D_1(k, s)$ is defined by (11).*

(i) *Then $k$-FWER $\leq \alpha$.*

(ii) *Moreover, for any stepup procedure with critical values of the form $\tilde{\alpha}_i = \alpha \alpha_i / D'$ for some constant $D'$ that satisfies $k$-FWER $\leq \alpha$, we have for each $i$ that $\alpha'_i \geq \tilde{\alpha}_i$.*

Before proceeding with the proof of Theorem 3.1, we recall the following lemma from [10], which generalizes an earlier result from [7].

LEMMA 3.2. *Suppose $\hat{p}_1, \ldots, \hat{p}_t$ are p-values in the sense that $P\{\hat{p}_i \leq u\} \leq u$ for all $i$ and $u$ in $(0, 1)$. Let their ordered values be $\hat{p}_{(1)} \leq \cdots \leq \hat{p}_{(t)}$. For some $m \leq t$, let*

$$0 = \beta_0 \leq \beta_1 \leq \beta_2 \leq \cdots \leq \beta_m \leq 1.$$

(i) *Then*

(12)
$$P\{\{\hat{p}_{(1)} \leq \beta_1\} \cup \{\hat{p}_{(2)} \leq \beta_2\} \cup \cdots \cup \{\hat{p}_{(m)} \leq \beta_m\}\}$$
$$\leq t \sum_{i=1}^{m}(\beta_i - \beta_{i-1})/i.$$

(ii) *As long as the right-hand side of* (12) *is* $\leq 1$, *the bound is sharp in the sense that there exists a joint distribution for the p-values for which the inequality is an equality.*

PROOF OF THEOREM 3.1. (i) Combining Lemmas 3.1 and 3.2, we have that

$$k - FWER \leq P\left\{\bigcup_{k \leq j \leq |I|} \{\hat{q}_{(j)} \leq \alpha'_{s-|I|+j}\}\right\}$$
$$\leq |I|\frac{\alpha'_{s-|I|+k}}{k} + |I| \sum_{k < j \leq |I|} \frac{\alpha'_{s-|I|+j} - \alpha'_{s-|I|+j-1}}{j}$$
$$= \frac{\alpha}{D_1(k,s)} S_1(k, s, |I|) \leq \alpha.$$

(ii) Consider the following joint distribution of $p$-values. Denote by $|I|^*$ the value of $|I|$ maximizing $S_1(k, s, |I|)$. Let the $p$-values of the $s - |I|^*$ false hypotheses be identically equal to 0 (or just $< \alpha'_1$) and let the $p$-values of the $|I|^*$ true hypotheses be constructed according to part (ii) of Lemma 3.2 so that

$$P\left\{\bigcup_{k \leq j \leq |I|^*} \{\hat{q}_{(j)} \leq \tilde{\alpha}_{s-|I|^*+j}\}\right\} = \frac{\alpha}{D'} S_1(k, s, |I|^*) = \frac{D_1}{D'}\alpha,$$

where the second equality uses the fact that $\tilde{\alpha}_i = \alpha \alpha_i / D'$. For such a joint distribution of $p$-values, the event of rejecting $\geq k$ true hypotheses is equivalent to rejecting $\geq s - |I|^* + k$ hypotheses in total. So,

$$k - FWER = P\left\{\bigcup_{k \leq j \leq |I|^*} \{\hat{q}_{(j)} \leq \tilde{\alpha}_{s-|I|^*+j}\}\right\} = \frac{D_1}{D'}\alpha.$$

Thus, to ensure control of the $k$-FWER, it must be the case that $D_1 \leq D'$. It follows that, for each $i$, $\alpha'_i \geq \tilde{\alpha}_i$. □

Theorem 3.1(ii) shows that it is not possible to increase all of the critical values by any amount without violating control of the $k$-FWER. In this sense, part (ii) of the theorem represents a sort of weak optimality result.

TABLE 1
*Values of $D_1(k,s)$ for k-FWER control with $\alpha_i$ given by* (13) *and* (19)

|   | $k=1$ | | $k=2$ | | $k=3$ | |
|---|---|---|---|---|---|---|
| $s$ | (13) | (19) | (13) | (19) | (13) | (19) |
| 10 | 2.11 | 3.92 | 2.03 | 2.57 | 1.90 | 2.10 |
| 25 | 2.13 | 7.99 | 2.16 | 4.72 | 2.15 | 3.60 |
| 50 | 2.13 | 14.52 | 2.16 | 8.10 | 2.17 | 5.91 |
| 100 | 2.13 | 27.32 | 2.16 | 14.63 | 2.17 | 10.33 |
| 250 | 2.13 | 65.25 | 2.16 | 33.77 | 2.17 | 23.22 |
| 500 | 2.13 | 128.08 | 2.16 | 65.34 | 2.17 | 44.36 |
| 1000 | 2.13 | 253.41 | 2.16 | 128.17 | 2.17 | 86.35 |
| 2000 | 2.13 | 503.75 | 2.16 | 253.51 | 2.17 | 170.01 |
| 5000 | 2.13 | 1254.20 | 2.16 | 628.96 | 2.17 | 420.46 |

Hommel and Hoffman [8] and Lehmann and Romano [10] propose using constants proportional to

(13) $$\alpha_i = \begin{cases} \dfrac{k}{s}, & \text{if } i \leq k, \\ \dfrac{k}{s+k-i}, & \text{if } i > k, \end{cases}$$

as part of a stepdown procedure to control the $k$-*FWER* and showed that such a procedure using critical values $\alpha\alpha_i$ controls the $k$-*FWER* at level $\alpha$ under no assumptions on the joint distribution of the $p$-values. We can apply Theorem 3.1 to this choice of $\alpha_i$ to construct a stepup procedure that also controls the $k$-*FWER* under no restrictions on the joint distribution of the $p$-values. Table 1 displays for several different values of $k$ and $s$ the normalizing constant $D_1(k,s)$ of Theorem 3.1. Table 1 shows that the constants must be approximately halved to ensure control of the $k$-*FWER*. For example, in the case $s=1000$ and $k=3$, the optimizing value of $|I|$ is 39, yielding $D_1(3, 1000) = 2.1707$.

For control of the *FWER*, Hochberg [5] proposed using the stepup procedure with critical values given by (13) with $k=1$. These same constants were used by Holm [6] to control the *FWER*, but as part of a stepdown procedure. Hochberg argued that his procedure controls the *FWER* assuming that the $p$-values are independent. Sarkar and Chang [15] have shown that Hochberg's procedure also controls the *FWER* for certain forms of positively dependent $p$-values. So, it follows from Remark 2.1 that under such assumptions on the joint distribution of the $p$-values Hochberg's procedure is more powerful than the one proposed by Holm. Holm's procedure, however, controls the *FWER* under no assumptions on the joint distribution of the $p$-values, whereas our results show that this is not true of Hochberg's procedure. However, we show that by dividing the constants by $D_1(1,s)$, control of the *FWER* is restored.

REMARK 3.1. The notion of control that we consider demands that $k\text{-}FWER \leq \alpha$ for *all* $P$. This is sometimes referred to as *strong* control of $k$-*FWER* in order to distinguish it from a weaker (and not particularly useful for multiple testing) notion of control known as *weak* control of the $k$-*FWER*, where it is only required that the $k\text{-}FWER \leq \alpha$ for all $P$ satisfying $|I| = |I(P)| = s$, that is, when all hypotheses are true. The distinction between weak and strong control generalizes in an obvious way to measures of error control other than the $k$-*FWER*. It is interesting to note that to guarantee even weak control of the $k$-*FWER*, the constants $\alpha\alpha_i$, where $\alpha_i$ is defined by (13), must be approximately halved (at least for large $s$). To see this, first note that when $|I| = s$, the $k$-*FWER* is equivalent to the probability of rejecting $\geq k$ hypotheses altogether; that is,

$$
(14) \qquad P\left\{\bigcup_{k \leq j \leq s} \{\hat{p}_{(j)} \leq \alpha\alpha_j\}\right\}.
$$

Using Lemma 3.2, we know there exists a joint distribution of the $p$-values for which (14) is equal to

$$
\begin{aligned}
(15) \quad & \alpha\left(1 + s \sum_{k < j \leq s} \frac{\alpha_j - \alpha_{j-1}}{j}\right) \\
& = \alpha\left(1 + k \sum_{k \leq i < s} \frac{s}{(s+k-i)i(i+1)}\right) \\
& = \alpha\left(1 + k \sum_{k \leq i < s} \frac{1}{i(i+1)} + k \sum_{k \leq i < s} \frac{i-k}{(s+k-i)i(i+1)}\right) \\
& = \alpha\left(2 - \frac{k}{s} + k \sum_{k \leq i < s} \frac{i-k}{(s+k-i)i(i+1)}\right).
\end{aligned}
$$

But, it is easy to see that, as $s \to \infty$, we have that (15) $\to 2\alpha$. It follows that, at least for large values of $s$, the constants $\alpha\alpha_i$ must be approximately halved to ensure weak control of the $k$-*FWER*. In fact, the expression (15) is strictly larger than the limiting value $2\alpha$, and so the constants must be divided by something slightly greater than 2. In order to guarantee strong control, the constants must be divided by something that is only slightly larger than 2. In the case $k = 1$, this value is 2.1314.

REMARK 3.2. More generally, suppose $|I|$ is not neccesarily $= s$ and denote the ordered values of the true $p$-values by $\hat{q}_{(1)} \leq \cdots \leq \hat{q}_{(|I|)}$. Then, following the argument given in the proof of Theorem 3.1(ii), we have that

$$k - FWER = P\left\{ \bigcup_{k \leq j \leq |I|} \{\hat{q}_{(j)} \leq \alpha\alpha_{s-|I|+j}\} \right\}. \tag{16}$$

Again, Lemma 3.2 asserts that there exists a joint distribution of true $p$-values for which (16) is equal to

$$\alpha\left(1 + |I| \sum_{k < j \leq |I|} \frac{\alpha_{s-|I|+j} - \alpha_{s-|I|+j-1}}{j}\right). \tag{17}$$

Note that

$$\alpha_{s-|I|+j} = \frac{k}{|I| + k - j},$$

so we may use the analysis of Remark 3.1 with the role of $s$ replaced by $|I|$ to conclude that (17) is equal to

$$\alpha\left(2 - \frac{k}{|I|}\right) + O\left(k\frac{\log |I|}{|I|}\right). \tag{18}$$

If $|I|$ is large, then it is sufficient to halve the constants $\alpha\alpha_i$ to control the $k$-$FWER$ approximately. The expression (18) implies further that if we index both $k$ and $|I|$ by the number of hypotheses $s$ and allow $s \to \infty$, then the stepup procedure with critical values $\alpha\alpha_i/2$ provides strong control of the $k$-$FWER$, provided that $k\frac{\log |I|}{|I|} \to 0$. Division by 2 can be thought of as the price to pay for using a stepup versus stepdown procedure (based on the same set of critical values). It is perhaps surprising that the value of 2 is independent of the choice of $k$.

Finner and Roters [3] compared stepup and stepdown procedures for control of the $FWER$ assuming that the $p$-values are exchangeable. Under the setup of their paper, their results suggest that stepup procedures are more powerful than stepdown procedures because one can use very nearly the same critical values for both procedures to control the $FWER$. However, in our comparisons, we assume nothing about the joint distribution of $p$-values and find that in such a setting the stepup procedure requires smaller critical values (by roughly a half) to provide control of the $FWER$, and more generally of the $k$-$FWER$.

We may also apply Theorem 3.1 to the sequence of constants given by

$$\alpha_i = \frac{i}{s}. \tag{19}$$

The normalizing constant $D_1(k, s)$ for this choice of $\alpha_i$ is also displayed in Table 1. In light of part (ii) of Theorem 3.1, we should not expect either of the sequences of critical values generated by applying Theorem 3.1 to (13) and (19) to be uniformly larger (and thus unambiguously more powerful) than the other. In order to illustrate this fact, we plot the two sequences of constants for the case in which $k = 2$, $s = 100$ and $\alpha = 0.05$. Panel (a) of Figure 1 displays the constants based on (13), whereas panel (b) displays the constants based on (19). Panel (c) depicts the ratio of the constants in panel (a) with the constants in panel (b). The dashed horizontal line in panel (c) is of height 1, allowing us to see graphically when the constants from panel (a) are greater than the constants from panel (b) and vice versa. We find that, for high and low values of $i$, the constants based on (13) are larger than the constants based on (19). For intermediate values of $i$, where the constants based on (13) are smaller than the constants based on (19), the differences between the constants are quite small in absolute terms, whereas, for other values of $i$, the differences between the constants are fairly substantial. This suggests that the procedure based on (13) may be preferable to the one based on (19).

**4. Control of the FDP.** The number $k$ of false rejections that one is willing to tolerate will often increase with the number of hypotheses rejected. This leads to consideration of not the number of false rejections (sometimes called false discoveries), but rather the proportion of false discoveries. Formally, let the *false discovery proportion* (*FDP*) be defined by

$$(20) \quad FDP = \begin{cases} \dfrac{\text{Number of false rejections}}{\text{Total number of rejections}}, & \text{if the denominator is } > 0, \\ 0, & \text{if there are no rejections.} \end{cases}$$

*FDP* is therefore the proportion of rejected hypotheses that are rejected erroneously. When none of the hypotheses is rejected, both numerator and denominator of that proportion are 0; since, in particular, there are no false rejections, the *FDP* is then defined to be 0. We now establish a general result that will aid us in constructing stepup procedures that control the *FDP* in the sense of (1). In what follows, we will sometimes use $m(j)$ as shorthand for $\lfloor \gamma j \rfloor + 1$, where $\lfloor x \rfloor$ is the greatest integer $\leq x$, and the notation $x \vee y$ in place of $\max\{x, y\}$.

LEMMA 4.1. *Consider testing $s$ null hypotheses, with $|I|$ of them true. Let $\hat{q}_{(1)} \leq \cdots \leq \hat{q}_{(|I|)}$ denote the ordered p-values corresponding to true hypotheses. Then the stepup procedure based on constants $\alpha_1 \leq \cdots \leq \alpha_s$ satisfies*

$$(21) \quad \begin{aligned} &P\{FDP > \gamma\} \\ &\leq P\left\{ \bigcup_{|I|-s+1 \leq k \leq |I|, |I| \geq m(s-|I|+k)} \{\hat{q}_{(k \vee m(s-|I|+k))} \leq \alpha_{s-|I|+k}\} \right\}. \end{aligned}$$

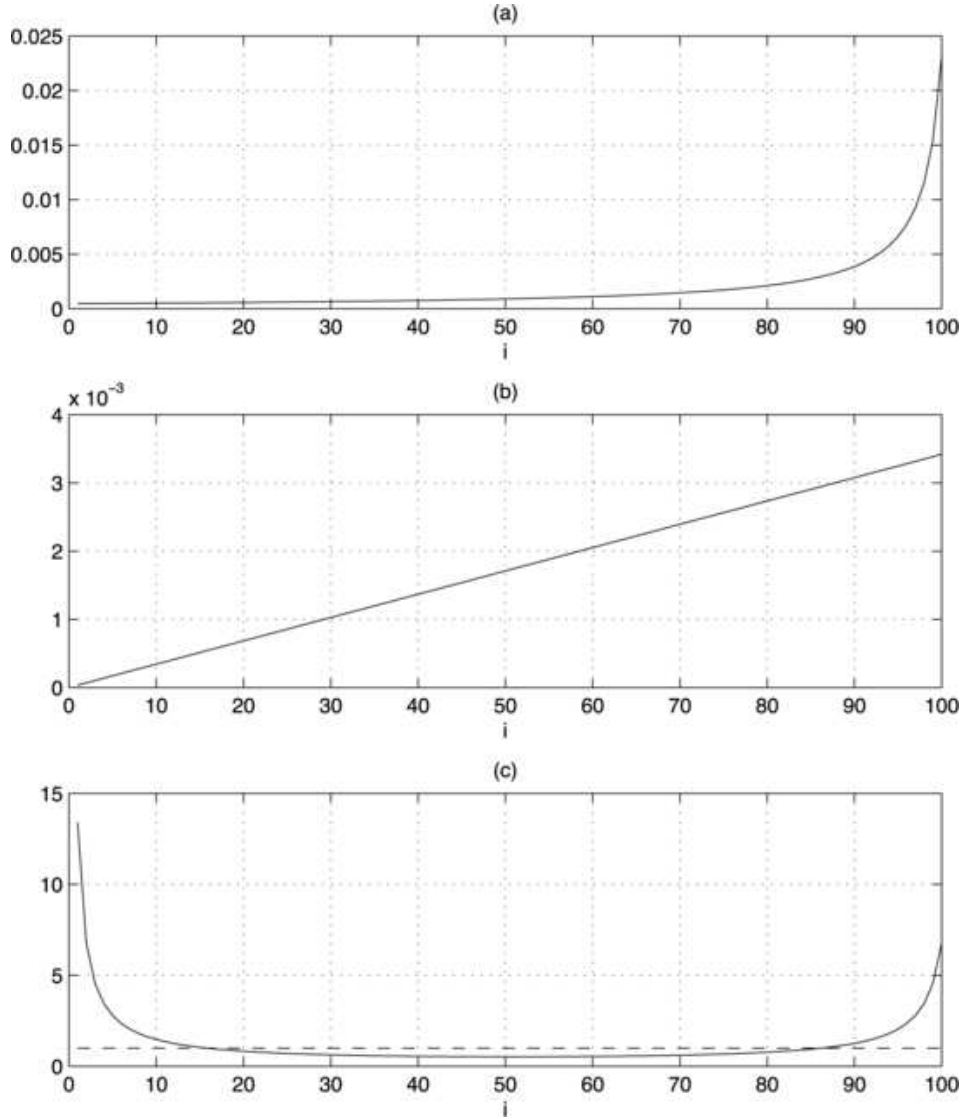

FIG. 1. *Stepup constants for k-FWER control with $k=2$, $s=100$ and $\alpha = 0.05$.*

PROOF. Let $A_j$ be the event in which exactly $j$ hypotheses are rejected by the stepup procedure; that is,

$$A_j = \{\hat{p}_{(s)} > \alpha_s, \ldots, \hat{p}_{(j+1)} > \alpha_j, \hat{p}_{(j)} \leq \alpha_j\}.$$

Let $T_j$ be the event in which at least $m(j)$ true hypotheses are rejected. Then

$$A_s \cap T_s \subseteq \{\hat{p}_{(s)} \leq \alpha_s\} \cap T_s \begin{cases} \subseteq \{\hat{q}_{(|I|)} \leq \alpha_s\}, & \text{if } |I| \geq m(s), \\ = \varnothing, & \text{otherwise.} \end{cases}$$

Likewise, for the event $A_{s-1} \cap T_{s-1}$, we have that

$$A_{s-1} \cap T_{s-1} \subseteq \{\hat{p}_{(s-1)} \leq \alpha_{s-1}\} \cap T_{s-1}$$

$$\begin{cases} \subseteq \{\hat{q}_{(|I|-1)} \leq \alpha_{s-1}\}, & \text{if } |I| - 1 \geq m(s-1), \\ \subseteq \{\hat{q}_{(|I|)} \leq \alpha_{s-1}\}, & \text{if } |I| = m(s-1), \\ = \varnothing, & \text{otherwise.} \end{cases}$$

It follows that

$$A_{s-1} \cap T_{s-1} \subseteq \{\hat{p}_{(s-1)} \leq \alpha_{s-1}\} \cap T_{s-1}$$

$$\begin{cases} \subseteq \{\hat{q}_{((|I|-1) \vee m(s-1))} \leq \alpha_{s-1}\}, & \text{if } |I| \geq m(s-1), \\ = \varnothing, & \text{otherwise.} \end{cases}$$

Following a similar line of reasoning, we have in general that

$$A_{s-j} \cap T_{s-j} \subseteq \{\hat{p}_{(s-j)} \leq \alpha_{s-j}\} \cap T_{s-j}$$

$$\begin{cases} \subseteq \{\hat{q}_{((|I|-j) \vee m(s-j))} \leq \alpha_{s-j}\}, & \text{if } |I| \geq m(s-j), \\ = \varnothing, & \text{otherwise.} \end{cases}$$

Thus

$$\{FDP > \gamma\} \subseteq \bigcup_{0 \leq j \leq s-1, |I| \geq m(s-j)} \{\hat{q}_{((|I|-j) \vee m(s-j))} \leq \alpha_{s-j}\}$$

$$= \bigcup_{|I|-s+1 \leq k \leq |I|, |I| \geq m(s-|I|+k)} \{\hat{q}_{(k \vee m(s-|I|+k))} \leq \alpha_{s-|I|+k}\},$$

from which the asserted claim follows. □

Given a sequence of constants $\alpha_1 \leq \cdots \leq \alpha_s$, we will now use Lemma 4.1 to construct a stepup procedure that satisfies (1). To this end, define

(22)
$$S_2 = S_2(\gamma, s, |I|)$$
$$= |I|\alpha_1 + |I| \sum_{|I|-s+1 < k \leq |I|, |I| \geq m(s-|I|+k)} \frac{\alpha_{s-|I|+k} - \alpha_{s-|I|+k-1}}{k \vee m(s-|I|+k)}$$

and let

(23) $$D_2 = D_2(\gamma, s) = \max_{1 \leq |I| \leq s} S_2(\gamma, s, |I|).$$

THEOREM 4.1. *Let $\alpha_1 \leq \cdots \leq \alpha_s$ be given. For testing $H_i : P \in \omega_i$, $i = 1, \ldots, s$, suppose $\hat{p}_i$ satisfies (3). Consider the stepup procedure with critical values $\alpha_i'' = \alpha \alpha_i / D_2(\gamma, s)$, where $D_2(\gamma, s)$ is defined by (23).*

(i) *Then $P\{FDP > \gamma\} \leq \alpha$; that is, (1) is satisfied.*

(ii) *Moreover, for any stepup procedure with critical values of the form $\tilde{\alpha}_i = \alpha \alpha_i / D'$ for some constant $D'$ that satisfies (1), we have for each $i$ that $\alpha_i'' \geq \tilde{\alpha}_i$.*

PROOF. (i) Combining Lemmas 4.1 and 3.2, we have that

$$P\{FDP > \gamma\} \leq P\left\{\bigcup_{|I|-s+1\leq k\leq |I|, |I|\geq m(s-|I|+k)} \{\hat{q}_{(k\vee m(s-|I|+k))} \leq \alpha''_{s-|I|+k}\}\right\}$$

$$\leq |I|\alpha_1'' + |I| \sum_{|I|-s+1 < k \leq |I|, |I| \geq m(s-|I|+k)} \frac{\alpha''_{s-|I|+k} - \alpha''_{s-|I|+k-1}}{k \vee m(s-|I|+k)}$$

$$= \frac{\alpha}{D_2(\gamma, s)} S_2(\gamma, s, |I|) \leq \alpha.$$

(ii) Consider the following joint distribution of $p$-values. Denote by $|I|^*$ the value of $|I|$ maximizing $S_2(\gamma, s, |I|)$. Let the distribution of the $p$-values corresponding to the $|I|^*$ true hypotheses be constructed according to part (ii) of Lemma 3.2 so that

$$P\left\{\bigcup_{|I|^*-s+1\leq k\leq |I|^*, |I|^*\geq m(s-|I|^*+k)} \{\hat{q}_{(k\vee m(s-|I|^*+k))} \leq \tilde{\alpha}_{s-|I|^*+k}\}\right\}$$

$$= \frac{\alpha}{D'} S_2(k, s, |I|^*) = \frac{D_2}{D'}\alpha,$$

where the second equality uses the fact that $\tilde{\alpha}_i = \alpha \alpha_i / D'$. We will now construct the joint distribution of the $p$-values corresponding to the $s - |I|^*$ false hypotheses conditional on the values of the true $p$-values so that $FDP > \gamma$ whenever (21) occurs. For the time being, suppose that $|I|^*$ is such that $|I|^* \geq m(s)$. Thus, (21) can be written more simply as

$$P\{FDP > \gamma\} \leq P\left\{\bigcup_{|I|^*-s+1\leq k\leq |I|^*} \{\hat{q}_{(k\vee m(s-|I|^*+k))} \leq \alpha''_{s-|I|^*+k}\}\right\}.$$

Define $k^*$ to be the smallest index $k > 0$ such that $k \geq m(s - |I|^* + k)$. Consider the event

$$(24) \quad \bigcup_{|I|^*\geq k\geq k^*} \{\hat{q}_{(k\vee m(s-|I|^*+k))} \leq \alpha_{s-|I|^*+k}\} = \bigcup_{|I|^*\geq k\geq k^*} \{\hat{q}_{(k)} \leq \alpha_{s-|I|^*+k}\}.$$

Whenever the event (24) occurs, let all false $p$-values be identically equal to 0. By assumption, $k \geq m(s - |I|^* + k)$ and $k > 0$, so note that whenever this event occurs, we have that $FDP > \gamma$.

Now suppose that the event (24) does not occur. Note that this rules out the possibility of any event of the form

$$\{\hat{q}_{(k \vee m(s-|I|^*+k))} \leq \alpha_{s-|I|^*+k}\}$$

for $k < k^*$ and $k \vee m(s - |I|^* + k) = k^*$. So, let $k^{**}$ be the largest $k < k^*$ such that $k \vee m(s - |I|^* + k) = k^* - 1$ and consider the event

(25) $$\{\hat{q}_{(k^*-1)} \leq \alpha_{s-|I|^*+k^{**}}\}.$$

Whenever the event (25) occurs but (24) does not, let $s - |I|^* + k^{**} - k^* + 1$ of the false $p$-values be identically equal to 0 and let the remaining $k^* - k^{**} - 1$ false $p$-values fall between $\alpha_{s-|I|^*+k^{**}}$ and $\alpha_{s-|I|^*+k^*}$. Again, by construction, whenever (24) does not occur but (25) does occur, we have $FDP > \gamma$.

We may continue arguing along these lines by replacing the role of $k^*$ with $k^{**}$ to construct a joint distribution of false $p$-values conditional on the true $p$-values such that, whenever (21) occurs, we have that $FDP > \gamma$. But we have assumed so far that $|I|^* \geq m(s)$. To generalize the argument to the case in which $|I|^* < m(s)$, note that the event (21) is always of the form

$$\bigcup_{1 \leq k \leq |I|^*} \{\hat{q}_{(k)} \leq \alpha_{l(k)}\}$$

for some strictly increasing sequence of positive integers $l(1) < \cdots < l(|I|^*)$. Thus, the smallest $l(|I|^*)$ can be is $|I|^*$. Let $(s - |I|^*) - (s - l(|I|^*)) = l(|I|^*) - |I|^* \geq 0$ of the false $p$-values be identically equal to 1. Since these hypotheses will always be accepted by the stepup procedure, we can restrict attention to the situation in which there are $s - l(|I|^*) + |I|^*$ hypotheses altogether, $s - l(|I|^*)$ of which are false. But for this situation, our assumption on the number of true hypotheses holds, so we may use the construction above to determine the distribution of remaining false $p$-values.

So, for such a joint distribution of $p$-values, we have that

$$P\{FDP > \gamma\} = P\left\{\bigcup_{|I|^*-s+1 \leq k \leq |I|^*, |I|^* \geq m(s-|I|^*+k)} \{\hat{q}_{(k \vee m(s-|I|^*+k))} \leq \tilde{\alpha}_{s-|I|^*+k}\}\right\}$$

$$= \frac{D_2}{D'}\alpha.$$

Thus, to ensure control the $FDP$, it must be the case that $D_2 \leq D'$. It follows that, for each $i$, $\alpha_i'' \geq \tilde{\alpha}_i$. □

Lehmann and Romano [10] develop a stepdown procedure that controls the $FDP$ in the sense of (1) by reasoning as follows. Denote by $F$ the number of false rejections. At step $i$, having rejected $i - 1$ hypotheses, we want to guarantee $F/i \leq \gamma$, that is, $F \leq \lfloor \gamma i \rfloor$. So, if $k = \lfloor \gamma i \rfloor + 1$, then $F \geq k$ should have probability no greater than $\alpha$; that is, we must control the number of false rejections to be $\leq k$. This leads them to consider using the stepdown constants (13) for control of the $k$-FWER with this particular choice of $k$ (which now depends on $i$). That is,

$$\alpha_i = \frac{\lfloor \gamma i \rfloor + 1}{s + \lfloor \gamma i \rfloor + 1 - i}. \tag{26}$$

Lehmann and Romano [10] provide two results that show that the stepdown procedure with critical values $\alpha \alpha_i$ with this choice of $\alpha_i$ satisfies (1). Unfortunately, some assumption on the joint dependence structure of the $p$-values is required. However, they show that if one considers a stepdown procedure with critical values $\alpha \alpha_i / C_{\lfloor \gamma s \rfloor + 1}$, where

$$C_j = \sum_{i=1}^{j} \frac{1}{i},$$

then the FDP is controlled in the sense of (1) without any assumptions on the dependence structure of the $p$-values.

Romano and Shaikh [12] show that this procedure is more conservative than necessary to control the FDP. Specifically, they show that the stepdown procedure with critical values obtained by replacing $C_{\lfloor \gamma s \rfloor + 1}$ with a smaller quantity $D_3(\gamma, s)$ also provides control of the FDP without any assumptions on the joint distribution of the $p$-values. This change leads to a considerable improvement, resulting in critical values typically 50 percent larger.

We can apply Theorem 4.1 to $\alpha_i$ defined by (26) to construct a stepup procedure that controls the FDP in the sense of (1). The normalizing constant $D_2(\gamma, s)$ is computed for several different values of $\gamma$ and $s$ in Table 2. The column labeled "$D_2$, (26)" refers to the value of $D_2$ when the constants (26) are used. For the purposes of comparison, we also display $D_3(\gamma, s)$. For large values of $s$, the normalizing constant $D_2(\gamma, s)$ is strictly smaller than $C_{\lfloor \gamma s \rfloor + 1}$, but it is always larger than $D_3(\gamma, s)$. Thus, it follows from Remark 2.1 that, for large values of $s$, the stepup procedure is more powerful than the stepdown procedure proposed by Lehmann and Romano [10], whereas a clear ranking of the procedure relative to the stepdown procedure proposed by Romano and Shaikh [12] is not possible.

As before with the $k$-FWER, we may also apply Theorem 4.1 to the sequence of constants defined by (19). The normalizing constant $D_2(\gamma, s)$ for this choice of $\alpha_i$ is also displayed in Table 2. Again, the optimality result stated in part (ii) of Theorem 4.1 suggests that we should not expect either of the sequences of critical values generated by applying Theorem 4.1 to (26) and (19) to be uniformly larger (and thus unambiguously more powerful) than the other. We plot the two sequences of constants for the special case in which $\gamma = 0.1$, $s = 100$ and $\alpha = 0.05$. Panel (a) of Figure 2 displays the constants based on (26), whereas panel (b) displays the constants based on (19). Panel (c) depicts the ratio of the constants in panel (a) with the constants in panel (b). The dashed horizontal line in panel (c) is of height 1, allowing us to see graphically when the constants from panel (a) are greater than the constants from panel (b) and vice versa. We find that, for high and low values of $i$, the constants based on (26) are larger than the constants based on (19). But, as with the comparison of the $k$-FWER controlling procedures, we find that the differences are comparatively small when the ones based on (26) are smaller than the constants based on (19) and fairly large otherwise. Thus, we believe the procedure based on (26) is likely to be preferred to the one based on (19).

Table 2
*Stepup constants for FDP control with $\alpha_i$ given by* (26) *and* (19)

| $s$ | $\gamma = 0.05$ | | | $\gamma = 0.1$ | | |
|---|---|---|---|---|---|---|
| | $D_2$, (26) | $D_2$, (19) | $D_3$, (26) | $D_2$, (26) | $D_2$, (19) | $D_3$, (26) |
| 10 | 2.11 | 3.91 | 1.00 | 2.11 | 3.91 | 1.00 |
| 25 | 2.40 | 7.99 | 1.43 | 2.68 | 7.78 | 1.50 |
| 50 | 2.70 | 14.12 | 1.50 | 2.99 | 10.96 | 1.75 |
| 100 | 2.96 | 20.32 | 1.73 | 3.37 | 15.09 | 2.04 |
| 250 | 3.41 | 31.04 | 2.12 | 3.93 | 21.21 | 2.52 |
| 500 | 3.80 | 40.33 | 2.50 | 4.39 | 26.33 | 2.95 |
| 1000 | 4.24 | 50.40 | 2.92 | 4.89 | 31.75 | 3.42 |
| 2000 | 4.72 | 61.05 | 3.38 | 5.41 | 37.37 | 3.92 |
| 5000 | 5.39 | 75.80 | 4.044 | 6.14 | 45.06 | 4.62 |

REMARK 4.1. Benjamini and Yekutieli [2] propose using the constants $\alpha \alpha_i / C_s$ for $\alpha_i$ given in (19) as part of a stepup procedure to control the *FDR* and show that such a procedure controls the *FDR* for all possible distributions of $p$-values. Since $FDR = FWER$ when $|I| = s$, we have that these critical values also control the *FWER* when $|I| = s$. But, the results in Table 1 show that these constants do not control the *FWER* in general since $D_1(1, s) > C_s$ for this choice of $\alpha_i$. More surprising, however, is that this observation continues to be true even if one assumes that the $p$-values are independent. To see this, consider the case in which $s$ is even and $|I| = s/2 + 1$. Suppose that all of the false $p$-values are identically equal to 0, and the true $p$-values, whose ordered values are denoted by $\hat{q}_{(1)} \leq \cdots \leq \hat{q}_{(|I|)}$, are each $\sim U(0, 1)$. Thus,

$$FWER = P\left\{ \bigcup_{1 \leq j \leq s/2+1} \left\{ \hat{q}_{(j)} \leq \frac{\alpha}{C_s} \alpha_{s/2-1+j} \right\} \right\}$$

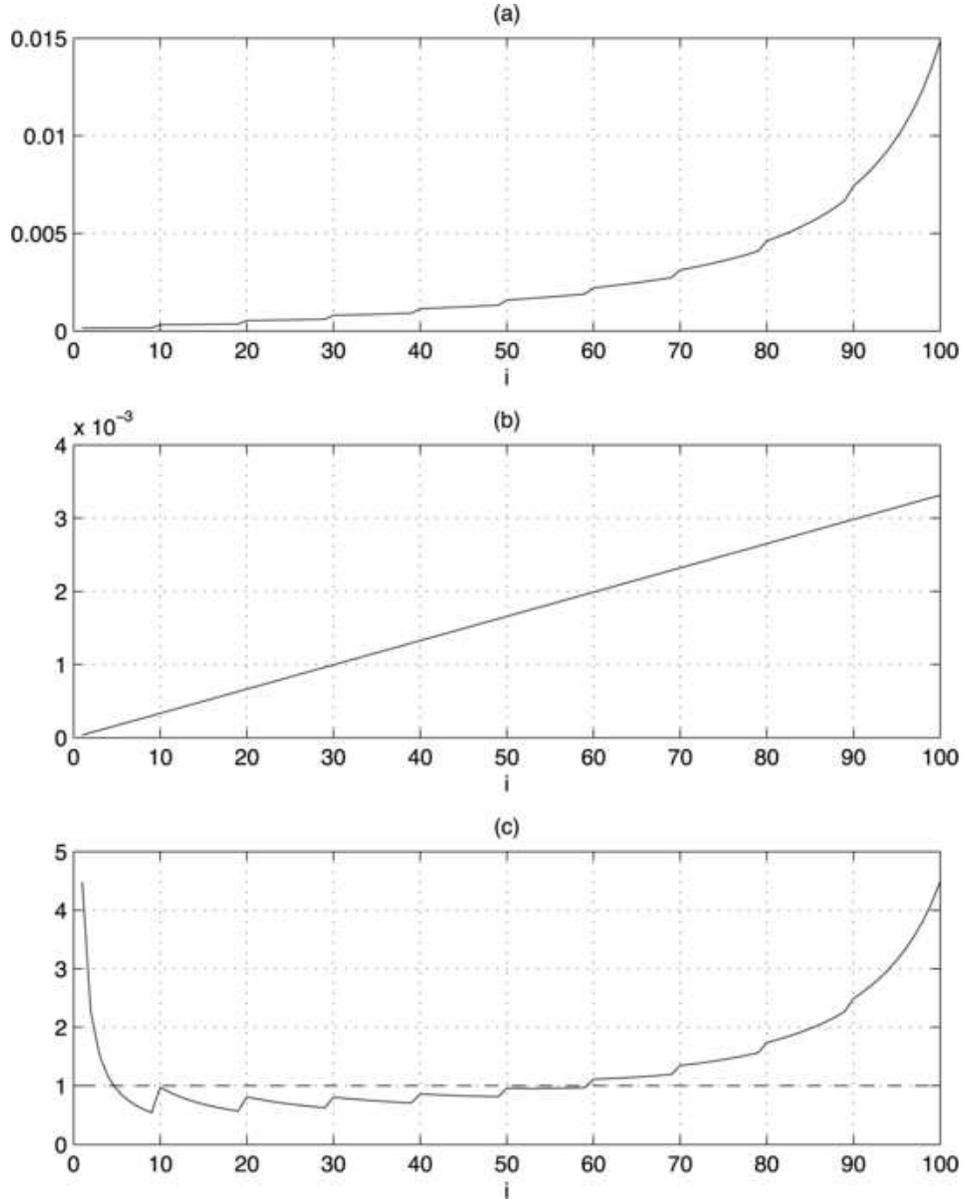

FIG. 2. *Stepup constants for FDP control with $\gamma = 0.1$, $s = 100$ and $\alpha = 0.05$.*

$$\geq P\left\{\hat{q}_{(1)} \leq \frac{\alpha}{C_s}\alpha_{s/2}\right\} = 1 - \left(1 - \frac{\alpha}{2C_s}\right)^{s/2+1} \to 1.$$

Thus, the stepup procedure with critical values $\alpha\alpha_i/C_s$ does not control the *FWER*, even under independence.

REMARK 4.2. In many situations the true individual $p$-values do not satisfy (3) exactly. However, suppose the $p$-values $\hat{p}_i^{(n)}$ are now indexed by $n$ (typically the sample size), and assume

$$\lim_{n\to\infty} P\{\hat{p}_i^{(n)} \leq u\} \leq u \qquad \text{for any } u, P \in \omega_i. \tag{27}$$

For example, if the $p$-values are determined by an asymptotic method such as the bootstrap, then it is typically the case that $\hat{p}_i^{(n)}$ converges in distribution to the uniform distribution on $(0,1)$ if $H_i$ is true. If we use a stepup procedure that controls the $FDP$ for nominal values of $\alpha$ and $\gamma$ whenever the $p$-values satisfy (3) exactly, then we can claim limiting control if we use a stepup procedure based on $p$-values which only satisfy (27). Specifically, we claim that asymptotic control holds; that is,

$$\limsup_{n\to\infty} P\{FDP(\hat{p}^{(n)}) > \gamma\} \leq \alpha, \tag{28}$$

where the event $\{FDP(\hat{p}^{(n)}) > \gamma\}$ that the $FDP$ is not controlled now shows the dependence on $n$ in that we are applying the procedure to the approximate vector of $p$-values $\hat{p}^{(n)} = (\hat{p}_1^{(n)}, \ldots, \hat{p}_s^{(n)})$. To see why, let $\hat{q}^{(n)} = (\hat{q}_1^{(n)}, \ldots, \hat{q}_{|I|}^{(n)})$ denote the $p$-values corresponding to the true hypotheses, with ordered values $\hat{q}_{(1)}^{(n)} \leq \cdots \leq q_{(|I|)}^{(n)}$. Then, by Lemma 4.1,

$$P\{FDP(\hat{p}^{(n)}) > \gamma\} \leq P\left\{\bigcup_k \hat{q}_{(k)}^{(n)} \leq \beta_k\right\}$$

for some nondecreasing $\beta_k$. We can write the right-hand side as $P\{\hat{q}^{(n)} \in C\}$, where $C$ is a closed set. (Note that the event that the $FDP$ is not controlled, viewed as a set in $s$-dimensional space, is not a closed set; there is no contradiction since the set $C$ corresponds to a larger set where the $FDP$ is not controlled.) But, $q^{(n)}$ is a tight sequence in $|I|$-dimensional Euclidean space (since it is supported on a fixed compact set, the $|I|$-fold product of $[0,1]$). So, taking any subsequence $\{n_j\}$, there exists a further subsequence $\{n_{j_l}\}$ along which $\hat{q}^{(n)}$ converges in distribution to a random vector $\hat{q} = (\hat{q}_1, \ldots, \hat{q}_{|I|})$ (which could depend on the subsubsequence). Moreover, the assumption (27) implies (3) holds for each $\hat{q}_i$. Let the ordered values of $\hat{q}$ be denoted $\hat{q}_{(1)} \leq \cdots \leq \hat{q}_{(|I|)}$. By the Portmanteau theorem, it follows that

$$\limsup_{n_{j_l}\to\infty} P\{\hat{q}^{(n_{j_l})} \in C\} \leq P\{q \in C\} = P\left\{\bigcup_k \hat{q}_{(k)} \leq \beta_k\right\}.$$

But, the right-hand side here is bounded above by $\alpha$ by Theorem 4.1. Since the bound $\alpha$ holds along any subsequence, the result is proved. A similar remark holds for control of the $k$-$FWER$ when using $p$-values that only satisfy (27) instead of (3).

TABLE 3
*Minimum and maximum values of ratios of Benjamini–Yekutieli constants and constants based on* (26) *and* (19) *when both are used to control the FDP*

| | $\gamma = 0.05$ | | | $\gamma = 0.1$ | | |
|---|---|---|---|---|---|---|
| $s$ | min (26) | max (26) | (19) | min (26) | max (26) | (19) |
| 10 | 9.25 | 27.76 | 14.96 | 4.63 | 13.88 | 7.48 |
| 25 | 4.71 | 31.86 | 9.55 | 2.33 | 14.26 | 4.91 |
| 50 | 2.75 | 33.40 | 6.37 | 1.99 | 15.05 | 4.10 |
| 100 | 2.25 | 35.02 | 5.11 | 1.86 | 15.41 | 3.44 |
| 250 | 2.03 | 35.80 | 3.93 | 1.75 | 15.53 | 2.88 |
| 500 | 1.95 | 35.72 | 3.368 | 1.68 | 15.46 | 2.58 |
| 1000 | 1.88 | 35.30 | 2.97 | 1.62 | 15.30 | 2.36 |
| 2000 | 1.81 | 34.67 | 2.68 | 1.58 | 15.10 | 2.19 |
| 5000 | 1.73 | 33.74 | 2.40 | 1.52 | 14.82 | 2.02 |

**5. Comparisons of *FDP* and *FDR* control.** In the previous section we have put forward two stepup procedures [one based on (26) and another based on (19)] that control the *FDP* in the sense of (1) under no assumptions on the dependence structure of the $p$-values. In this section we will use the crude inequalities given in (2) to compare these two *FDP*-controlling procedures with the *FDR*-controlling stepup procedure of Benjamini and Yekutieli [2].

From the second inequality in (2), control of the *FDR* at level $\gamma\alpha$ implies control of the *FDP* in the sense of (1). Since the Benjamini and Yekutieli stepup procedure with constants $\alpha i/(sC_s)$ controls the *FDR* at level $\alpha$, the constants given by

$$\alpha'_i = \frac{\gamma\alpha i}{sC_s} \tag{29}$$

control the *FDP* under no assumptions on the joint distribution of the $p$-values. We will first compare these critical values with the critical values of the form $\alpha\alpha_i/D_2(\gamma, s)$ derived by applying Theorem 4.1 to $\alpha_i$ defined by (26). Note that the ratio of the critical values $\alpha\alpha_i/D_2(\gamma, s)$ to $\alpha'_i$ is only a function of $\gamma$ and $s$. Table 3 displays for several different values of $\gamma$ and $s$ the minimum and maximum values of this ratio. For all values of $\gamma$ and $s$ in the table, the minimum value of the ratio $> 1$. In fact, the value of $\alpha\alpha_i/D_2(\gamma, s)$ is often at least twice as large as the corresponding value of $\alpha'_i$. The procedure based on the constants (26) is therefore unambiguously more powerful than the procedure based on the constants (29). By examining the maximum value of the ratio, we see that the value of $\alpha\alpha_i/D_2(\gamma, s)$ may be more than 15 times as large as the corresponding value of $\alpha'_i$.

We may replace the critical values based on (26) with those based on (19) and perform the same comparison. In this case, the ratio of $\alpha\alpha_i/D_2(\gamma, s)$ to $\alpha'_i$ is simply $C_s/(\gamma D_2(\gamma, s))$ and does not depend on $i$. Table 3 also displays the value of this ratio for several values of $\gamma$ and $s$. We find that the critical values $\alpha \alpha_i / D_2(\gamma, s)$ are always at least twice as large as the critical values $\alpha'_i$; thus, as before, the procedure based on the constants (19) is more powerful than the procedure based on the constants (29).

It is also possible to utilize the *FDP*-controlling constants to control the *FDR*, by application of (2). Now using (29) results in larger critical values than those resulting from application of Theorem 4.1. Detailed numerical comparisons are available from the authors. These results, though based on the crude inequalities in (2), suggest that it is perhaps worthwhile to consider the sort of control desired when choosing critical values. Indeed, the previous comparisons are somewhat unfair in that the *FDR*-controlling procedures were not designed to control the *FDP*, and vice versa.

However, we consider one final comparison in which the *FDP*-controlling constants are utilized to control the median of the *FDP* at level $\gamma$ by setting $\alpha = 1/2$. We may compare these critical values with the Benjamini–Yekutieli critical values given by $\alpha''_i = \gamma i / s C_s$, which control the *FDR* at level $\gamma$. First, we consider the constants based on (26). Table 4 displays the minimum and maximum values of the ratio of these critical values to the critical values $\alpha''_i$ for several different values of $\gamma$ and $s$. We find that, for moderate values of $s$, the critical values based on (26) are uniformly larger than the critical values $\alpha''_i$, but, for large values of $s$, the critical values $\alpha''_i$ are larger for some values of $i$. To examine whether these differences are of any practical significance, we plot in Figure 3 the two sequences of constants for the case in which $s = 1000$ and $\gamma = 0.1$. Panel (a) displays the critical values based on (26), whereas panel (b) displays the critical values $\alpha''_i$. Panel (c) displays the ratio of the constants in panel (a) with the constants in panel (b). The dashed horizontal line in panel (c) is of height 1. It is clear that, except for some small values of $i$, the constants of panel (a) are often dramatically larger than the constants of panel (b). More importantly, at such values of $i$, the differences between the two sequences of critical values are quite small. Thus, for most practical purposes, the stepup procedure based on the constants in panel (a) seems preferable to the one based on the constants in panel (b).

We now consider the same comparison with the critical values based on (26) replaced by the critical values based on (19). For this choice of $\alpha_i$, the value of the ratio of the constants derived from Theorem 4.1 to the constants $\alpha''_i$ no longer depends on $i$. Table 4 displays the values of this ratio for several values of $\gamma$ and $s$. Here we find that the critical values based on (19) used to control the median of the *FDP* are always uniformly larger, and therefore more powerful, than the *FDR*-controlling critical values $\alpha''_i$, though, for large values of $s$, the two sequences of critical values are nearly indistinguishable.

TABLE 4
*Minimum and maximum values of ratios of Benjamini–Yekutieli constants when used to control the FDR and constants based on* (26) *and* (19) *when used to control the median of the FDP*

| $s$ | $\gamma = 0.05$ | | | $\gamma = 0.1$ | | |
|---|---|---|---|---|---|---|
|  | min (26) | max (26) | (19) | min (26) | max (26) | (19) |
| 10 | 4.63 | 13.88 | 7.48 | 2.31 | 6.94 | 3.74 |
| 25 | 2.36 | 15.93 | 4.78 | 1.17 | 7.13 | 2.45 |
| 50 | 1.37 | 16.70 | 3.18 | 1.00 | 7.52 | 2.05 |
| 100 | 1.12 | 17.51 | 2.55 | 0.93 | 7.71 | 1.72 |
| 250 | 1.02 | 17.90 | 1.97 | 0.88 | 7.77 | 1.44 |
| 500 | 0.98 | 17.86 | 1.68 | 0.84 | 7.73 | 1.29 |
| 1000 | 0.94 | 17.65 | 1.49 | 0.81 | 7.65 | 1.18 |
| 2000 | 0.91 | 17.34 | 1.34 | 0.79 | 7.55 | 1.09 |
| 5000 | 0.87 | 16.87 | 1.20 | 0.76 | 7.411 | 1.01 |

## 6. Empirical applications.

EXAMPLE 6.1 (*Benjamini–Hochberg application*). We revisit the study of treatments for myocardial infarction analyzed in Benjamini and Hochberg [1], Section 3.2. For the 15 reported $p$-values, the Benjamini–Hochberg *FDR* controlling procedure at level 0.05 rejects four hypotheses. But, this procedure does not work for all possible joint distributions of $p$-values. The more generally applicable Benjamini–Yekutieli procedure rejects only three hypotheses. In contrast, our procedure for controlling the median of the *FDP* at level 0.05 still rejects four hypotheses.

EXAMPLE 6.2 (*Comparing strategies to a benchmark*). The problem considered is to determine which, if any, of several financial strategies outperforms a given benchmark. The data set is similar to that in Romano and Wolf [14]. We consider all $s = 210$ hedge funds in the Center for International Securities and Derivatives Markets (CISDM) database that have a complete return history from 01/1994 to 12/2003. All returns are net of management and incentive fees. The benchmark is the risk free rate of return. Performance is measured monthly, so each fund has a return history of 120 values. It is well known that returns of hedge funds exhibit nontrivial serial correlations and the distribution of (Studentized) differences in log returns between a particular strategy and a benchmark must take into account such dependence. Individual or marginal $p$-values were calculated according to the Studentized circular block bootstrap, as reviewed in [14].

For $k$-*FWER* control at level $\alpha = 0.1$ our stepup procedure rejects 10, 20 or 24 hypotheses according to $k = 1$, $k = 5$ or $k = 10$. For control of the *FDR*, the Benjamini–Yekutieli procedure rejects 0, 16 or 23 according to whether $\gamma = 0.01$, $\gamma = 0.05$ or $\gamma = 0.1$, respectively. For control of the median of the FDP, our procedure rejects 20, 22 or 24 hypotheses according to the same values of $\gamma = 0.01$, $\gamma = 0.05$ or $\gamma = 0.1$.

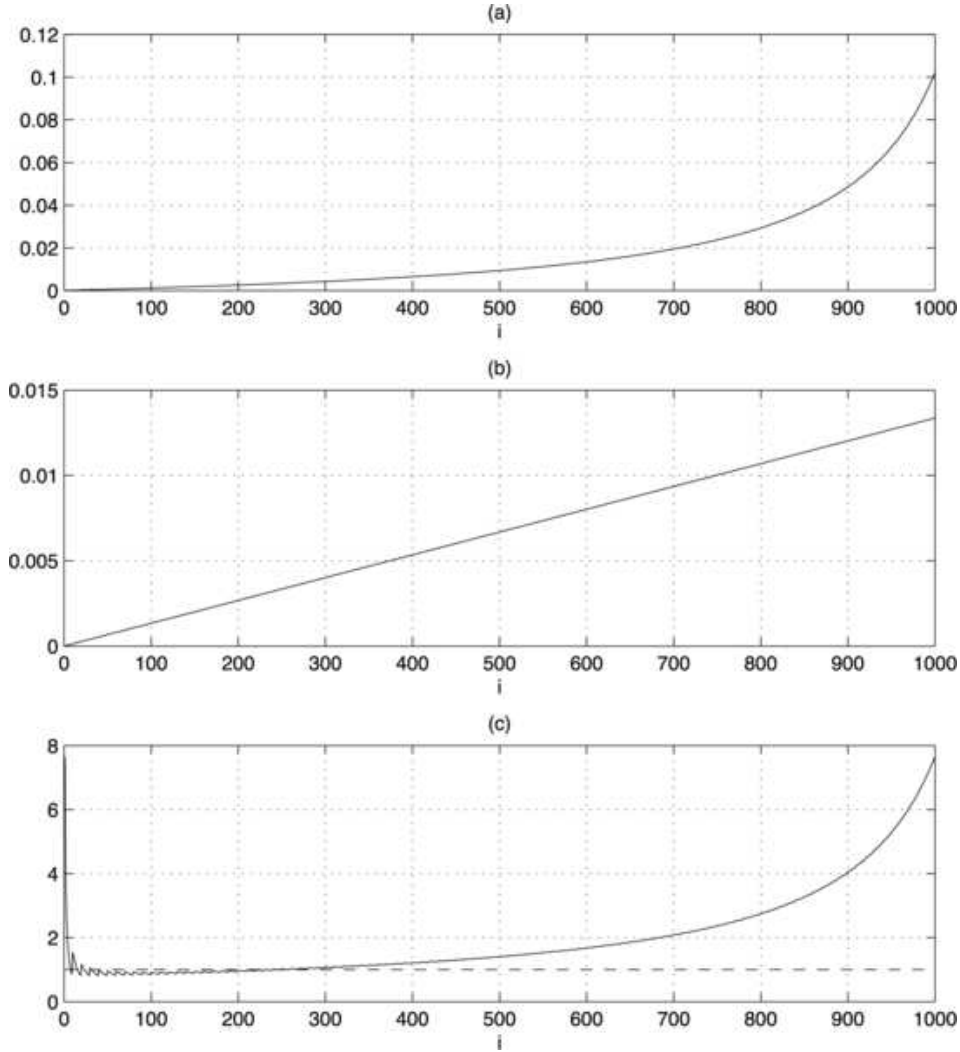

Fig. 3. *Stepup constants with median of FDP $\leq \gamma$ and FDR $\leq \gamma$ for $s = 1000$ and $\gamma = 0.1$.*

**7. Conclusion.** In this article we have described stepup procedures for testing multiple hypotheses that control either the $k$-FWER or the FDP without any restrictions on the joint distribution of the $p$-values. For each of these two measures of error control, we have also shown that the procedures constructed using our results satisfy a sort of weak optimality in that the critical values cannot all be made larger without violating the measure of error control. Our results have also revealed that control of the $k$-*FWER* or *FDP* using a stepup procedure assuming nothing about the joint distribution of $p$-values requires smaller critical values than a stepdown procedure satisfying the same measure of error control. Finally, we have compared two *FDP*-controlling procedures obtained using our results with the stepup procedure for control of the *FDR* of Benjamini and Yekutieli [2], which is also valid under no assumptions on the joint distribution of the $p$-values. These comparisons suggest that it is indeed important to consider the sort of error control desired when constructing multiple testing procedures.

**Acknowledgment.** Thanks to Michael Wolf for computation of the $p$-values in Example 6.2.

## REFERENCES

[1] BENJAMINI, Y. and HOCHBERG, Y. (1995). Controlling the false discovery rate: A practical and powerful approach to multiple testing. *J. Roy. Statist. Soc. Ser. B* **57** 289–300. MR1325392

[2] BENJAMINI, Y. and YEKUTIELI, D. (2001). The control of the false discovery rate in multiple testing under dependence. *Ann. Statist.* **29** 1165–1188. MR1869245

[3] FINNER, H. and ROTERS, M. (1998). Asymptotic comparison of step-down and step-up multiple test procedures based on exchangeable test statistics. *Ann. Statist.* **26** 505–524. MR1626043

[4] GENOVESE, C. and WASSERMAN, L. (2004). A stochastic process approach to false discovery control. *Ann. Statist.* **32** 1035–1061. MR2065197

[5] HOCHBERG, Y. (1988). A sharper Bonferroni procedure for multiple tests of significance. *Biometrika* **75** 800–802. MR0995126

[6] HOLM, S. (1979). A simple sequentially rejective multiple test procedure. *Scand. J. Statist.* **6** 65–70. MR0538597

[7] HOMMEL, G. (1983). Tests of the overall hypothesis for arbitrary dependence structures. *Biometrical J.* **25** 423–430. MR0735888

[8] HOMMEL, G. and HOFFMAN, T. (1987). Controlled uncertainty. In *Multiple Hypothesis Testing* (P. Bauer, G. Hommel and E. Sonnemann, eds.) 154–161. Springer, Heidelberg.

[9] KORN, E., TROENDLE, J., MCSHANE, L. and SIMON, R. (2004). Controlling the number of false discoveries: Application to high-dimensional genomic data. *J. Statist. Plann. Inference* **124** 379–398. MR2080371

[10] LEHMANN, E. L. and ROMANO, J. P. (2005). Generalizations of the familywise error rate. *Ann. Statist.* **33** 1138–1154. MR2195631

[11] LEHMANN, E. L. and ROMANO, J. P. (2005). *Testing Statistical Hypotheses*, 3rd ed. Springer, New York. MR2135927

[12] ROMANO, J. P. and SHAIKH, A. M. (2004). On control of the false discovery proportion. Technical Report No. 2004-31, Dept. Statistics, Stanford Univ.

[13] ROMANO, J. P. and WOLF, M. (2005a). Control of generalized error rates in multiple testing. Technical Report 2005–2012. Dept. Statistics, Stanford Univ.

[14] ROMANO, J. P. and WOLF, M. (2005). Stepwise multiple testing as formalized data snooping. *Econometrica* **73** 1237–1282. MR2149247
[15] SARKAR, S. and CHANG, C. (1997). The Simes method for multiple hypothesis testing with positively dependent test statistics. *J. Amer. Statist. Assoc.* **92** 1601–1608. MR1615269
[16] VAN DER LAAN, M., DUDOIT, S. and POLLARD, K. (2004). Augmentation procedures for control of the generalized family-wise error rate and tail probabilities for the proportion of false positives. *Stat. Appl. Genet. Mol. Biol.* **3** Article 15. MR2101464

DEPARTMENT OF STATISTICS
STANFORD UNIVERSITY
STANFORD, CALIFORNIA 94305-4065
USA
E-MAIL: romano@stanford.edu

DEPARTMENT OF ECONOMICS
STANFORD UNIVERSITY
STANFORD, CALIFORNIA 94305-6072
USA
E-MAIL: ashaikh@stanford.edu